\documentclass{ifacconf}  
\usepackage{IEEEtrantools}
\usepackage{natbib}
\usepackage[utf8]{inputenc} % set input encoding (not needed with XeLaTeX)
\usepackage{mathrsfs}
\usepackage{amsfonts}
\usepackage{amssymb}
\usepackage{amsmath}
\usepackage{graphicx}
\usepackage{commath}
\usepackage{tikz}
\usepackage{tkz-graph}
\usetikzlibrary{arrows}
\usepackage{pgfplots}
\usetikzlibrary{external}
\newtheorem{theorem}{Theorem}

\newtheorem{remark}{Remark}
\newtheorem{objective}{Objective}

\DeclareMathOperator*{\diag}{diag}

\DeclareMathOperator*{\sgn}{sgn}

\newlength\fheight
\newlength\fwidth

\begin{document}
\begin{frontmatter}

\title{Distributed Voltage and Current Control of Multi-Terminal High-Voltage Direct Current Transmission Systems
\thanksref{footnoteinfo}} % Title, preferably not more than 10 words.

\thanks[footnoteinfo]{This work was supported in part by the European Commission by the Hycon2 project, the Swedish Research Council (VR) and the Knut and Alice Wallenberg Foundation. 
The $2^{\text{nd}}$ and $6^{\text{th}}$ authors are supported by ELEKTRA. 
The $3^{\text{rd}}$ author is also affiliated with the Centre for Autonomous Systems at KTH.
 Corresponding author: Martin Andreasson, e-mail: mandreas@kth.se.}

\author[First]{Martin Andreasson} 
\author[Second]{Mohammad Nazari} 
\author[First]{Dimos V. Dimarogonas}
\author[First]{Henrik Sandberg}
\author[First]{Karl H. Johansson}
\author[Second]{Mehrdad Ghandhari} 

\address[First]{ACCESS Linnaeus Centre, School of Electrical Engineering,
KTH Royal Institute of Technology, Sweden.\\
(e-mail: {\tt \small \{mandreas, dimos, hsan, kallej\}@kth.se}.) }

\address[Second]{Electric Power Systems, School of Electrical Engineering, KTH Royal Institute of Technology, Stockholm, Sweden.\\
(e-mail: {\tt \small \{nazarim, mehrdad\}@kth.se}.) }

\begin{abstract}   
High-voltage direct current (HVDC) is a commonly used technology for long-distance power transmission, due to its low resistive losses and low costs. 
In this paper, a novel distributed controller for multi-terminal HVDC (MTDC) systems is proposed. Under certain conditions on the controller gains, it is shown to stabilize the MTDC system. The controller is shown to always keep the voltages close to the nominal voltage, while assuring that the injected power is shared fairly among the converters. The theoretical results are validated by simulations, where the affect of communication time-delays is also studied. 
\end{abstract}

\end{frontmatter}

%%%%%%%%%%%%%%%%%%%%%%%%%%%%%%%%%%%%%%%%
%%%%%%%%%%%%%%%%%%%%%%%%%%%%%%%%%%%%%%%%
\section{Introduction}

Transmitting power over long distances is one of the greatest challenges in today's power transmission systems. Increased distances between power generation and consumption is a driving factor behind long-distance power transmission. High-voltage direct current (HVDC) is a commonly used technology for long-distance power transmission, due to its low resistive losses and lower costs compared to AC transmission systems. 
Off-shore wind farms also typically require HVDC power transmission, as the need for reactive current limits the maximum transmission capacity of AC power transmission lines.

 With increased HVDC line constructions, future HVDC transmission systems are likely to consist of multiple terminals, to be able to connect several AC systems. 
 Voltage source converters make it possible to build HVDC systems with multiple terminals, referred to as multi-terminal HVDC (MTDC) systems in the literature. 

Maintaining an adequate DC voltage is one of the most important control problem for MTDC transmission systems. If the DC voltage deviates too far from the nominal operational voltage, equipment could be damaged \citep{18}. 

Different voltage control methods for MTDC systems have been proposed in the literature. Among them, the voltage margin method (VMM) and the voltage droop method (VDM) are the most well-known methods \citep{8}. These control methods change the injected active power from the alternating current (AC) systems into the DC grid to maintain active power balance in the DC grid and as a consequence, control the DC voltage. A decreasing DC voltage requires increased injected currents through the converters in order to restore the voltage.

VDM is designed so that all or more than one converter participate to control the DC voltage \citep{14}. All participant terminals change their injected active power to control the DC voltage. A higher slope of the voltage characteristic means that a terminal will inject less power given a certain change in the DC voltage.

VMM on the other hand, is designed so that one terminal is responsible to control the DC voltage, while the other terminals keep their injected active power constant. The terminal controlling the DC voltage is referred to as the slack terminal. When the slack terminal is no longer able to supply or extract the power necessary to maintain its DC bus voltage above a certain voltage margin, a new terminal will operate as the slack terminal \citep{8}. 

 A promising alternative approach to control MTDC networks is to use various distributed voltage controllers instead of VDM or VMM controllers. 
Distributed control has been successfully applied to both primary and secondary frequency control of AC transmission systems  \citep{Andreasson2012, Andreasson2013_ecc, simpson2012synchronization}. Recently, various distributed controllers have been applied also to voltage control of MTDC transmission systems \citep{nazari1}, including distributed secondary frequency control of asynchronous AC transmission systems \citep{dai2010impact}. In this paper, we propose a novel distributed voltage controller for MTDC transmission systems, which possesses the property of power sharing.

%\begin{figure*}[ht!]
%    \centering
%    $
%\begin{array}{cc}
%    {\includegraphics[height=2.0in,width=2.75in]{VDM.pdf}}
%& 
%    {\includegraphics[height=2.0in,width=2.75in]{vmm.pdf}}
%    \\
%\text{(a)} & \text{(b)}
%\end{array}
%$
%    \caption{(a) VDM characteristic for a two-terminal HVDC system. (b) VMM characteristic for a two-terminal HVDC. If the active power is positive, the active power direction is from AC grid to the DC grid and the converter operates at rectifier mode; otherwise, it operates at inverter mode. } \label{fig:VC}
%\end{figure*}

 This remainder of this paper is organized as follows. In Section \ref{sec:prel}, the mathematical notation is defined. In Section \ref{sec:model}, the system model and the control objectives are defined. In Section \ref{sec:MTDC-control}, a voltage droop controller is presented and analysed. Subsequently, a distributed averaging controller is presented, and its stability and steady-state properties. In Section \ref{sec:simulations}, simulations of the distributed controller on a four-terminal MTDC test system are provided, before ending with a discussion and concluding remarks in Section~\ref{sec:discussion}.

\section{Notation}
\label{sec:prel}
Let $\mathcal{G}$ be an undirected graph. Denote by $\mathcal{V}=\{ 1,\hdots, n \}$ the vertex set of $\mathcal{G}$, and by $\mathcal{E}=\{ 1,\hdots, m \}$ the edge set of $\mathcal{G}$. Let $\mathcal{N}_i$ be the set of neighboring vertices to $i \in \mathcal{V}$.
In this paper we will only consider static and connected graphs. For the application of control of MTDC power transmission systems, this is a reasonable assumption as long as there are no power line failures.
Denote by $\mathcal{B} = \mathcal{B(G)}$ the vertex-edge adjacency matrix of $\mathcal{G}$, and let $\mathcal{\mathcal{L}_W}=\mathcal{B}W\mathcal{B}^T$ be the weighted Laplacian matrix of $\mathcal{G}$, with edge-weights given by the  elements of the diagonal matrix $W$.
Let $\mathbb{C}^-$ denote the open left half complex plane, and $\bar{\mathbb{C}}^-$ its closure. We denote by $c_{n\times m}$ a vector or matrix of dimension $n\times m$ whose elements are all equal to $c$. $I_{n}$ denotes the identity matrix of dimension $n$. For simplicity, we will often drop the notion of time dependence of variables, i.e., $x(t)$ will be denoted $x$.

%%%%%%%%%%%%%%%%%%%%%%%%%%%%%%%%%%%%%%%%
%%%%%%%%%%%%%%%%%%%%%%%%%%%%%%%%%%%%%%%%
\section{Model and problem setup}
\label{sec:model}
\begin{figure}
\begin{center}
\begin{tikzpicture}
\GraphInit[vstyle=Normal]
% Vertices
\Vertex[x=0, y=0] {1}
\Vertex[x=2, y=0] {2}
\Vertex[x=0, y=-1.5] {3}
\Vertex[x=2, y=-1.5] {4}
% Edgesfr
\Edge[label=$e_1$](1)(2)
\Edge[label=$e_4$](3)(4)
\Edge[label=$e_2$](1)(3)
\Edge[label=$e_3$](2)(4)
\end{tikzpicture}
\end{center}
\caption{Example of a graph topology of a MTDC system.}
\label{fig:graph}
\end{figure}
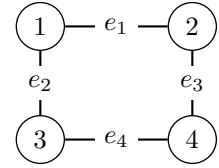

Consider an MTDC transmission system consisting of $n$ converters, denoted $1, \dots, n$, see Figure \ref{fig:graph} for an example of an MTDC topology. The converters are assumed to be connected by $m$ HVDC transmission lines. The dynamics of converter $i$ is assumed to be given by
\begin{align}
\begin{aligned}
C_i \dot{V}_i &= -\sum_{j\in \mathcal{N}_i} I_{ij} + I_i^{\text{inj}} + u_i \\
&= -\sum_{j\in \mathcal{N}_i} \frac{1}{R_{ij}}(V_i -V_j) + I_i^{\text{inj}} + u_i,
\end{aligned}
\label{eq:hvdc_scalar}
\end{align}
where $V_i$ is the voltage of converter $i$, $C_i$ is its capacity, $I_i^{\text{inj}}$ is the nominal injected current, which is assumed to be unknown but constant over time, and $u_i$ is the controlled injected current.  The constant $R_{ij}$ denotes the resistance of the transmission line connecting the converters $i$ and $j$. Equation \eqref{eq:hvdc_scalar} may be written in vector-form as
\begin{align}
\begin{aligned}
\dot{V} &= -C\mathcal{L}_R V +CI^{\text{inj}} + Cu,
\end{aligned}
\label{eq:hvdc_vector}
\end{align}
where $V=[V_1, \dots, V_n]^T$, $C=\diag([C_1^{-1}, \dots, C_n^{-1}])$, $I^{\text{inj}} = [I^{\text{inj}}_1, \dots, I^{\text{inj}}_n]^T$, $u=[u_1, \dots, u_n]^T$ and $\mathcal{L}_R$ is the weighted Laplacian matrix of the graph representing the transmission lines, denoted $\mathcal{G}_R$, whose edge-weights are given by the conductances $\frac{1}{R_{ij}}$. The control objectives considered in this paper are twofold.

\begin{objective}
\label{obj:1}
The voltages of the converters, $V_i$, should converge to a value close to the nominal voltage $V^{\text{nom}}$, after a disturbance has occurred. The nominal voltage $V^{\text{nom}}$ is assumed to be identical for all converters. It is however clear that it is not possible to have $\lim_{t\rightarrow \infty}V_i(t) = V^{\text{nom}}$ for all $i \in \mathcal{V}$, since this would imply that the currents between all converters are zero.
\end{objective}

\begin{objective}
\label{obj:2}
The injected currents should converge to a value which is proportional to an a priori known parameter, i.e.
\begin{align*}
\lim_{t\rightarrow \infty} u(t) = K^u 1_{n\times 1},
\end{align*}
for some diagonal matrix $K^u$, whose elements are positive.
 The second objective is often referred to as \emph{power sharing}, in the sense that the ratios between the injected currents of the converters is always the same at stationarity. Since $P_{i}^{\text{inj}}=V_iI_{i}^{\text{inj}}$, and since the relative voltage differences of the converters are very small, the injected power can be well approximated as being proportional to the injected current. 
\end{objective}

%%%%%%%%%%%%%%%%%%%%%%%%%%%%%%%%%%%%%%%%
%%%%%%%%%%%%%%%%%%%%%%%%%%%%%%%%%%%%%%%%
\section{MTDC control}
\label{sec:MTDC-control}

%%%%%%%%%%%%%%%%%%%%%%%%%%%%%%%%%%%%%%%%
\subsection{Voltage droop control}
\label{subsec:model_droop}
In this section the voltage droop method (VDM) will be studied, as well as some of its limitations. VDM is a simple decentralized proportional controller taking the form
\begin{align}
u_i &= K^P_i(V^{\text{nom}}-V_i),
\label{eq:droop}
\end{align}
where $V^{\text{nom}}$ is the nominal DC voltage. Alternatively, the controller \eqref{eq:droop} can be written in vector form as
\begin{align}
u &= K^P(V^{\text{nom}}1_{n\times 1}-V),
\label{eq:droop_vector}
\end{align}
where $K^P=\diag([K^P_1, \dots, K^P_n])$.
The decentralized structure of the voltage droop controller is often advantageous for control of HVDC converters, as the time constant of the voltage dynamics is typically smaller than the communication delays between the converters. The DC voltage regulation is typically carried out by all converters. However, the VDM possesses some drawbacks. Firstly, the voltages of the converters don't converge to a value close to the nominal voltages in general. 
Secondly, the controlled injected currents do not converge to a certain ration, i.e., power sharing, as shown in the following theorem. % \ref{th:droop_stability}. 
%(WE MUST SHOW IT. BUT THIS IS WHAT WE SHOW WITH THEOREM 1, RIGHT?)

\begin{theorem}
\label{th:droop_stability}
Consider an MTDC transmission system described by \eqref{eq:hvdc_scalar}, where the control input $u_i$ is given by \eqref{eq:droop} and the injected currents $I_i^{\text{inj}}$ are constant. The closed-loop system is stable for any $K^P>0$, in the sense that the voltages $V$ converge to some constant value. However $\lim_{t\rightarrow \infty} V_i(t) \ne V^{\text{nom}}$ in general. Furthermore, the controlled injected currents satisfy $\lim_{t\rightarrow \infty} \sum_{i=1}^n (u^i+I_i^{\text{inj}}) = 0$. However $\lim_{t\rightarrow  \infty } u(t) \ne -(\sum_{i=1}^n I_i^{\text{inj}})/(\sum_{i=1}^n K^P_i) K^u 1_{n\times 1}$ in general, for any diagonal $K^u$ with positive elements.
\end{theorem}

\begin{pf}
The closed loop dynamics of \eqref{eq:hvdc_vector} with $u$ given by \eqref{eq:droop} are
\begin{align}
\begin{aligned}
\dot{V} &= -C\mathcal{L}_R V + CK^P(V^{\text{nom}}1_{n\times 1}-V) + CI^{\text{inj}} \\
&= \underbrace{-C(\mathcal{L}_R + K^P)}_{\triangleq A}V + CK^PV^{\text{nom}}1_{n\times 1} + CI^{\text{inj}}.
\end{aligned}
\label{eq:hvdc_closed_loop_droop_vector}
\end{align}
Clearly the dynamics \eqref{eq:hvdc_closed_loop_droop_vector} are stable if and only if $A$ as defined above is Hurwitz. Consider the characteristic polynomial of A:
\begin{align*}
\begin{aligned}
0&= \det(sI_n-A) = \det\left(sI_n + C(\mathcal{L}_R + K^P)\right) \Leftrightarrow \\
0&= \det \underbrace{\left(sC^{-1} + (\mathcal{L}_R + K^P)\right)}_{\triangleq Q(s)}.
\end{aligned}
\end{align*}
The equation $0=\det(Q(s))$ has a solution for a given $s$ only if $0=x^TQ(s)x$ has a solution for some $\norm{x}=1$. This gives
\begin{align*}
0&=s\underbrace{x^TC^{-1}x}_{a_1} + \underbrace{x^T(\mathcal{L}_R + K^P)x}_{a_0}.
\end{align*}
Clearly $a_0, a_1>0$, which implies that the above equation has all its solutions $s\in \mathbb{C}^-$ by the Routh-Hurwitz stability criterion. This implies that the solutions of $0=\det(Q(s))$ satisfy $s\in \mathbb{C}^-$, and thus that $A$ is Hurwitz.

Now consider the equilibrium of \eqref{eq:hvdc_closed_loop_droop_vector}:
\begin{align}
0&= {-C(\mathcal{L}_R + K^P)}V + CK^PV^{\text{nom}}1_{n\times 1} + CI^{\text{inj}}.
\label{eq:hvdc_closed_loop_droop_vector_equilibrium}
\end{align}
Since $K^P>0$ by assumption $(\mathcal{L}_R + K^P)$ is invertible, which implies
\begin{align}
V=(\mathcal{L}_R + K^P)^{-1}\left( K^PV^{\text{nom}}1_{n\times 1} + I^{\text{inj}} \right),
\label{eq:hvdc_closed_loop_droop_vector_equilibrium_voltage}
\end{align}
which does not equal to $V^{\text{nom}}1_{n\times 1}$ in general. It is also easily seen that
$$ u \ne (\sum_{i=1}^n I_i^{\text{inj}})/(\sum_{i=1}^n K^P_i) K^P 1_{n\times 1}$$ in general. Premultiplying \eqref{eq:hvdc_closed_loop_droop_vector_equilibrium} with $1_{1\times n} C^{-1}$ yields
\begin{align*}
0 &= 1_{1\times n} K^P(V^{\text{nom}}1_{n\times 1} -V) + I^{\text{inj}} = \sum_{i=1}^n (u_i + I_i^{\text{inj}}) \qed
\end{align*} 
\end{pf}
Generally when tuning the proportional gains $K^P$, there is a trade-off between having the voltages converge to the nominal voltage, and having power sharing between the converters. Having low gains $K^P$ will result in better power sharing properties, but the voltages will be far from the reference value. On the other hand, having high gains  $K^P$ will ensure that the voltages converge close to the nominal voltage, at the expense of the power sharing properties. This rule of thumb is formalized in the following theorem.

\begin{theorem}
\label{th:droop_power_sharing}
Consider an MTDC network described by \eqref{eq:hvdc_scalar}, where the control input $u_i$ is given by \eqref{eq:droop} with positive gains $K^P_i$, and constant injected currents $I_i^{\text{inj}}$. The DC voltages satisfy
\begin{align*}
\lim_{K_i^P\rightarrow \infty \; \forall i=1, \dots, n} \lim_{t\rightarrow \infty} V(t) &= V^{\text{nom}}1_{n\times 1} \\
\lim_{K_i^P\rightarrow 0 \; \forall i=1, \dots, n} \lim_{t\rightarrow \infty} V(t) &= \sgn\left(\sum_{i=1}^n I_i^{\text{inj}}\right)\infty 1_{n\times 1},
\end{align*}
while the controlled injected currents satisfy
\begin{align*}
&\lim_{K_i^P\rightarrow \infty \; \forall i=1, \dots, n} \lim_{t\rightarrow \infty} u(t) = -I^{\text{inj}} \\
&\lim_{K_i^P\rightarrow 0 \; \forall i=1, \dots, n} \lim_{t\rightarrow \infty} u(t) = \\  &-\left(\sum_{i=1}^n I_i^{\text{inj}}\right)/\left(\sum_{i=1}^n K^P_i\right) K^P 1_{n\times 1},
\end{align*}
\end{theorem}

\begin{pf}
Let us first consider the case when $K_I^P\rightarrow \infty \; \forall i=1, \dots, n$. In the equilibrium of \eqref{eq:hvdc_closed_loop_droop_vector}, the voltages satisfy by \eqref{eq:hvdc_closed_loop_droop_vector_equilibrium_voltage}:
\begin{align*}
& \lim_{K_i^P\rightarrow \infty \; \forall i=1, \dots, n}  V
\\ &= \lim_{K_i^P\rightarrow \infty \; \forall i=1, \dots, n} (\mathcal{L}_R + K^P)^{-1}\left( K^PV^{\text{nom}}1_{n\times 1} + I^{\text{inj}} \right) \\
&= \lim_{K_i^P\rightarrow \infty \; \forall i=1, \dots, n} ( K^P)^{-1}\left( K^PV^{\text{nom}}1_{n\times 1} + I^{\text{inj}} \right) \\
&= V^{\text{nom}}1_{n\times 1}.
\end{align*}
By inserting the above expression for the voltages, the controlled injected currents are given by
\begin{align*}
&\lim_{K_i^P\rightarrow \infty \; \forall i=1, \dots, n} u = \lim_{K_i^P\rightarrow \infty \; \forall i=1, \dots, n} K^P\left( V^{\text{nom}}1_{n\times 1} -V   \right) \\
&= \lim_{K_i^P\rightarrow \infty \; \forall i=1, \dots, n} K^P \left( - (K^P)^{-1}I^{\text{inj}}  \right){=}{-}I^{\text{inj}}.
\end{align*}

Now consider the case when $\lim_{K_i^P\rightarrow 0 \; \forall i=1, \dots, n}$. Since $(\mathcal{L}_R + K^P)$ is real and symmetric, any vector in $\mathbb{R}^n$ can be expressed as a linear combination of its eigenvectors. Denote by $(v_i, \lambda_i)$ the eigenvector and eigenvalue pair $i$ of $(\mathcal{L}_R + K^P)$. Write
\begin{align}
\left( K^PV^{\text{nom}}1_{n\times 1} + I^{\text{inj}} \right) &= \sum_{i=1}^n a_i v_i,
\label{eq:eigendecomposition_1}
\end{align}
where $a_i, i=1, \dots, n$ are real constants. The equilibrium of \eqref{eq:hvdc_closed_loop_droop_vector} implies that the voltages satisfy
\begin{align*}
&\lim_{K_i^P\rightarrow 0 \; \forall i=1, \dots, n}  V \\
&= \lim_{K_i^P\rightarrow 0 \; \forall i=1, \dots, n} (\mathcal{L}_R + K^P)^{-1}\left( K^PV^{\text{nom}}1_{n\times 1} + I^{\text{inj}} \right) \\
&= \lim_{K_i^P\rightarrow 0 \; \forall i=1, \dots, n} (\mathcal{L}_R + K^P)^{-1}\sum_{i=1}^n a_i v_i \\&=\lim_{K_i^P\rightarrow 0 \; \forall i=1, \dots, n} \sum_{i=1}^n \frac{a_i}{\lambda_i} v_i = \frac{a_1}{\lambda_1} v_1,
\end{align*}
where $\lambda_1$ is the smallest eigenvalue of $(\mathcal{L}_R + K^P)$, which clearly satisfies $\lambda_1 \rightarrow 0^+$ as ${K_i^P\rightarrow 0 \; \forall i=1, \dots, n}$. Hence the last equality in the above equation holds. By letting ${K_i^P\rightarrow 0 \; \forall i=1, \dots, n}$ and premultiplying   \eqref{eq:eigendecomposition_1} with $v_1^T=1/n1_{n\times 1}$, we obtain $a_1=(\frac{1}{n} \sum_{i=1}^n I_i^{\text{inj}})$ since the eigenvectors of $(\mathcal{L}_R + K^P)$ form an orthonormal basis of $\mathbb{R}^n$. Thus $\lim_{K_i^P\rightarrow 0 \; \forall i=1, \dots, n} \lim_{t\rightarrow \infty} V(t) = \sgn\left(\sum_{i=1}^n I_i^{\text{inj}}\right)\infty 1_{n\times 1}$. Finally the controlled injected currents are given by
\begin{align*}
&\lim_{K_i^P\rightarrow 0 \; \forall i=1, \dots, n} u = \lim_{K_i^P\rightarrow 0 \; \forall i=1, \dots, n} K^P(V^{\text{nom}}1_{n\times 1} - V) \\
&= \lim_{K_i^P\rightarrow 0 \; \forall i=1, \dots, n} K^P\big(V^{\text{nom}}1_{n\times 1} - \frac{a_1}{\lambda_1} 1_{n\times 1}\big) \\
&= -\frac{a_1}{\lambda_1} K^P 1_{n\times 1}.
\end{align*}
By premultiplying \eqref{eq:hvdc_closed_loop_droop_vector_equilibrium} with $1_{1\times n} C^{-1}$ we obtain
\begin{align*}
1_{1\times n} K^P(V^{\text{nom}}1_{n\times 1} - V) = - 1_{1\times n} I^{\text{inj}},
\end{align*}
which implies that
\begin{align*}
\frac{a_1}{\lambda_1} &= \frac{1_{1\times n} I^{\text{inj}}}{1_{1\times n}  K^P 1_{n\times 1} } = (\sum_{i=1}^n I_i^{\text{inj}})/(\sum_{i=1}^n K^P_i),
\end{align*}
which gives the desired expression for $u$. \hfill  $\qed$
\end{pf}

%%%%%%%%%%%%%%%%%%%%%%%%%%%%%%%%%%%%%%%%
\subsection{Distributed averaging control}
In this section we propose a distributed controller for MTDC transmission systems which allows for communication between the converters. The proposed controller takes inspiration from the control algorithms given by \cite{Andreasson2013_ecc, Andreasson2013_PI} and by \cite{nazari1}, and is given by
\begin{align}
\label{eq:distributed_voltage_control}
\begin{aligned}
u_i &= K_i^P(\hat{V}_i -V_i) \\
\dot{\hat{V}}_i &= K^V_i(V^{\text{nom}}{-}V_i)-\gamma \sum_{j\in \mathcal{N}_i} c_{ij} \left( (\hat{V}_i -V_i){-}(\hat{V}_j -V_j) \right),
\end{aligned}
\end{align}
where $\gamma>0$ is a constant, and
\begin{align*}
K^V_i &= \left\{ \begin{array}{ll}
1 & \text{if } i=1 \\
0 & \text{otherwise.}
\end{array} \right.
\end{align*}
This controller can be understood as a fast proportional control loop (consisting of the first line), and a slower integral control loop (consisting of the second line). The internal controller variables $\hat{V}_i$ can be understood as reference values for the proportional control loops, regulated by the integral control loop. 
Converter $i=1$, without loss of generality, acts as voltage regulator. The first line of \eqref{eq:distributed_voltage_control} ensures that the controlled injected currents are quickly adjusted after a change in the voltage. $c_{ij}=c_{ji}>0$ is a constant, and $\mathcal{N}_i$ denotes the set of converters which can communicate with converter $i$. The communication graph is assumed to be undirected, i.e., $j\in \mathcal{N}_i$ implies $i\in \mathcal{N}_j$. The second line ensures that the voltage is restored at converter $1$ by integral action, and that the controlled injected currents are proportional to the proportional gains $K^P_i$ at stationarity. In vector-form, \eqref{eq:distributed_voltage_control} can be written as
\begin{align}
\label{eq:distributed_voltage_control_vector}
\begin{aligned}
u &= K^P (\hat{V} -V) \\
\dot{\hat{V}} &= K^V (V^{\text{nom}}1_{n\times 1} -V) -\gamma \mathcal{L}_c (\hat{V}-V) ,
\end{aligned}
\end{align}
where $K^P$ is defined as before, $K^V=\diag([K^V_1, 0, \dots, 0])$, and $\mathcal{L}_C$ is the weighted Laplacian matrix of the graph representing the communication topology, denoted $\mathcal{G}_c$, whose edge-weights are given by $c_{ij}$, and which is assumed to be connected.
 The following theorem shows that the proposed controller \eqref{eq:distributed_voltage_control} has the desirable properties which the droop controller \eqref{eq:droop} is lacking, and gives sufficient conditions for which controller parameters result in a stable closed loop system.

\begin{theorem}
\label{th:distributed_voltage_control}
Consider an MTDC network described by \eqref{eq:hvdc_scalar}, where the control input $u_i$ is given by \eqref{eq:distributed_voltage_control} and the injected currents $I^{\text{inj}}$ are constant.
The closed loop system is stable if
\begin{align}
&\begin{aligned}
&\frac{1}{2}\lambda_{\min} \left( (K^P)^{-1}\mathcal{L}_R + \mathcal{L}_R(K^P)^{-1} \right) + 1 +\\
& \frac{\gamma}{2}\lambda_{\min} \left( \mathcal{L}_C (K^P)^{-1} C^{-1} + C^{-1}  (K^P)^{-1} \mathcal{L}_C \right)  >0
\end{aligned} \label{eq:stability_cond_distributed_1} \\
& \lambda_{\min} \left(  \mathcal{L}_C (K^P)^{-1} \mathcal{L}_R + \mathcal{L}_R (K^P)^{-1} \mathcal{L}_C \right) \ge 0. \label{eq:stability_cond_distributed_2}
\end{align}
Furthermore 
$$\lim_{t\rightarrow  \infty } u(t) = -(\sum_{i=1}^n I_i^{\text{inj}})/(\sum_{i=1}^n K^P_i) K^P 1_{n\times 1},$$ 
and  $\lim_{t\rightarrow \infty} V_1(t) = V^{\text{nom}}$. This implies that the controlled injected currents satisfy Objective \ref{obj:2}, with $K^u = (\sum_{i=1}^n I_i^{\text{inj}})/(\sum_{i=1}^n K^P_i) K^P$. The remaining voltages satisfy $\lim_{t\rightarrow \infty} |V_i(t)-V^{\text{nom}}|\le 2I^{\text{max}}\sum_{i=2}^n \frac{1}{\lambda_i}$, where $I^{\text{max}}=\max_i |I^{\text{tot}}|$ and $I^{\text{tot}}= \lim_{t\rightarrow \infty} I^{\text{inj}} +u(t)$. Here $\lambda_i$ denotes the $i$'th eigenvalue of $\mathcal{L}_R$. 
\end{theorem}

\begin{remark}
\label{rem:1}
There always exists a sufficiently large $K^P$, and sufficiently small $\gamma$, such that the condition \eqref{eq:stability_cond_distributed_1} is fulfilled.
\end{remark}

\begin{remark}
\label{rem:2}
A sufficient condition for when \eqref{eq:stability_cond_distributed_2} is fulfilled, is that $\mathcal{L}_c=k_2\mathcal{L}_k$, $k_2 \in \mathbb{R}^+$ i.e., the topology of the communication network is identical to the topology of the power transmission lines, up to a positive scaling factor.
\end{remark}

\begin{pf}
The closed loop dynamics of \eqref{eq:hvdc_vector} with the controlled injected currents $u$ given by \eqref{eq:distributed_voltage_control_vector} are given by
\begin{align}
\begin{bmatrix}
\dot{\hat{V}} \\ \dot{V}
\end{bmatrix}
&{=}
\underbrace{
\begin{bmatrix}
-\gamma \mathcal{L}_C & \gamma \mathcal{L}_C -K^V \\
CK^P & -C(\mathcal{L}_R +K^P)
\end{bmatrix}}_{\triangleq A}
\begin{bmatrix}
{\hat{V}} \\ {V}
\end{bmatrix}
{+}
\begin{bmatrix}
K^V V^{\text{nom}}1_{n\times 1}  \\
CI^{\text{inj}}
\end{bmatrix}.
\label{eq:closed_loop_distributed_vector}
\end{align}
The characteristic equation of $A$ is given by
\begin{align*}
0 &= \det(sI_{2n} -A) = \left| \begin{matrix}
sI_n +\gamma \mathcal{L}_C & -\gamma \mathcal{L}_C + K^V \\
-CK^P & sI_n + C(\mathcal{L}_R +K^P)
\end{matrix} \right| \\
&= \frac{|CK^P|}{|sI_n+\gamma \mathcal{L}_c|} 
\left| \begin{matrix}
sI_n +\gamma \mathcal{L}_C & -\gamma \mathcal{L}_C + K^V \\
-sI_n-\gamma\mathcal{L}_c & \begin{matrix}
(sI_n+\gamma\mathcal{L}_c)(K^P)^{-1}C^{-1}\cdot \\ (sI_n + C(\mathcal{L}_R +K^P))
\end{matrix} 
\end{matrix} \right| \\
&= |CK^P| |(sI_n+\gamma\mathcal{L}_c)(K^P)^{-1}C^{-1} (sI_n + C(\mathcal{L}_R +K^P)) \\ 
&-\gamma \mathcal{L}_C + K^V | = \\
&= |CK^P| \left( (\gamma \mathcal{L}_c(K^P)^{-1}\mathcal{L}_R +K^V) + s((K^P)^{-1}\mathcal{L}_R +I_n \right. \\
&+\left. \gamma\mathcal{L}_C(K^P)^{-1}C^{-1}) +s^2((K^P)^{-1}C^{-1})  \right)
\\
 &\triangleq |CK^P| \det(Q(s)).
\end{align*}
This assumes that $|sI_n+\gamma \mathcal{L}_c|\ne 0$, however $|sI_n+\gamma \mathcal{L}_c|= 0$ implies $s =0$ or $s\in \mathbb{C}^-$. However, since $A$ is full rank, this still implies that all solutions satisfy $s\in \mathbb{C}^-$. 
Now, the above equation has a solution only if $x^TQ(s)x=0$ for some $x:\; \norm{x}=1$. This condition gives the following equation
\begin{align*}
0&=  \underbrace{x^T(\gamma \mathcal{L}_C (K^P)^{-1} \mathcal{L}_R + K^V)x}_{a_0}  \\
& \;\;\;\;  + s\underbrace{x^T((K^P)^{-1}\mathcal{L}_R + I_n + \gamma \mathcal{L}_C(K^P)^{-1}C^{-1} )x}_{a_1} \\ 
& \;\;\;\; + s^2 \underbrace{x^T((K^P)^{-1}C^{-1})x}_{a_2},
\end{align*}
which by the Routh-Hurwitz stability criterion has all solutions $s\in \mathbb{C}^-$ if and only if $a_i>0$ for $i=0,1,2$. 

Clearly, $a_2>0$, since $((K^P)^{-1}C^{-1})$ is diagonal with positive elements. It is easily verified that $a_1>0$ if
\begin{align*}
&\frac{1}{2}\lambda_{\min} \left( (K^P)^{-1}\mathcal{L}_R + \mathcal{L}_R(K^P)^{-1} \right) \\
+& \frac{\gamma}{2}\lambda_{\min} \left( \mathcal{L}_C (K^P)^{-1} C^{-1} + C^{-1}  (K^P)^{-1} \mathcal{L}_C \right) + 1 >0.
\end{align*}
Finally, clearly $x^T( \mathcal{L}_C (K^P)^{-1} \mathcal{L}_R)x \ge 0$ for any $x:\; \norm{x}=1$ if and only if
\begin{align*}
\frac 12 \lambda_{\min} \left(  \mathcal{L}_C (K^P)^{-1} \mathcal{L}_R + \mathcal{L}_R (K^P)^{-1} \mathcal{L}_C \right) &\ge 0.
\end{align*}
Since the graphs corresponding to $\mathcal{L}_R$ and $\mathcal{L}_C$ are both assumed to be connected, the only $x$ for which $x^T( \mathcal{L}_C (K^P)^{-1} \mathcal{L}_R)x = 0$ is $x=\frac{1}{\sqrt{n}}[1, \dots, 1]^T$. Given this $x=\frac{1}{\sqrt{n}}[1, \dots, 1]^T$, $x^T K^Vx=\frac{1}{n} K^V_1>0$. Thus, $a_0>0$ given that the above inequality holds. Thus, under assumptions \eqref{eq:stability_cond_distributed_1}--\eqref{eq:stability_cond_distributed_2}, $A$ is Hurwitz, and thus the closed loop system is stable.

Now consider the equilibrium of \eqref{eq:closed_loop_distributed_vector}. Premultiplying the first $n$ rows with $1_{1\times n}$ yields $0=1_{1\times n} K^V(V^{\text{nom}}1_{n\times 1}-V) = K^V_1(V^{\text{nom}}-V_1)$. Inserting this back to the first $n$ rows of \eqref{eq:closed_loop_distributed_vector} yields $0=\mathcal{L}_C (V-\hat{V})$, implying that $(V-\hat{V}) = k1_{n\times 1}$. Inserting this in \eqref{eq:distributed_voltage_control_vector} gives $u=K^P(V-\hat{V}) = kK^P1_{n\times 1}$.
To obtain a bound on the remaining voltages, we consider again the equilibrium of  \eqref{eq:closed_loop_distributed_vector}. The last $n$ rows of the equilibrium of \eqref{eq:closed_loop_distributed_vector} give
\begin{align}
\label{eq:equilibrium_I_tot}
\mathcal{L}_R V &= K^P(\hat{V}-V) + I^{\text{inj}} = I^{\text{tot}}.
\end{align}
Let
\begin{align*}
V=\sum_{i=1}^n a_i w_i,
\end{align*}
where $w_i$ is the $i$'th eigenvector of $\mathcal{L}_R$ with the corresponding eigenvalue $\lambda_i$. Since $\mathcal{L}_R$ is symmetric, the eigenvectors $\{w_i\}_{i=1}^n$ can be chosen so that they form an orthonormal basis of $\mathbb{R}^n$. Using the eigendecomposition of $V$ above, we obtain the following equation from \eqref{eq:equilibrium_I_tot}:
\begin{align}
\label{eq:equilibrium_I_tot_eigen}
\mathcal{L}_R V &= \mathcal{L}_R \sum_{i=1}^n a_i w_i =  \sum_{i=1}^n a_i \lambda_i w_i = I^{\text{tot}}.
\end{align}
By premultiplying \eqref{eq:equilibrium_I_tot_eigen} with $w_k$ for $k=1, \dots , n$, we obtain:
\begin{align*}
a_k\lambda_k = w_k^T I^{\text{tot}},
\end{align*}
due to orthonormality of $\{w_i\}_{i=1}^n$. Hence, for $i=2, \dots, n$ we get
\begin{align*}
a_k = \frac{w_k^T I^{\text{tot}}}{\lambda_k}.
\end{align*}
The constant $a_1$ is however not determined by \eqref{eq:equilibrium_I_tot_eigen}, since $\lambda_1=0$. Denote $\Delta V=\sum_{i=2}^n a_i w_i$. Since $w_1=\frac{1}{\sqrt{n}}1_{n\times 1}$, $V_i-V_j=\Delta V_i-\Delta V_j$ for any $i,j \in \mathcal{V}$. Thus, the following bound is easily obtained:
\begin{align*}
&|V_i-V_j| = |\Delta V_i-\Delta V_j| \le 2 \max_i |\Delta V_i| = 2 \norm{\Delta V}_\infty \\
 &\le 2 \norm{\Delta V}_2 = 2 \norm{\sum_{i=2}^n a_i w_i}_2 \le 2 \sum_{i=2}^n |a_i|= 2 \sum_{i=2}^n \left|\frac{w_i^T I^{\text{tot}}}{\lambda_i}\right| \\
 & \le 2 I^{\text{max}} \sum_{i=2}^n \frac{1}{\lambda_i},
\end{align*}
where we have used the fact that $\norm{w_i}_2=1$ for all $i=1, \dots, n$, and $\norm{x}_\infty \le \norm{x}_2$ for any $x\in \mathbb{R}^n$. Since the upper bound on $|V_i-V_j|$ is valid for any $i,j \in \mathcal{V}$, it is in particular valid for $j=1$. Recalling that for the equilibrium $V_1=V^{\text{nom}}$, the desired inequality is obtained.
Finally, setting $\dot{V}=0_{n\times 1}$ in \eqref{eq:hvdc_vector} and premultiplying with $1_{1\times n} C^{-1}$ gives
$
0=1_{1\times n}I^{\text{inj}} + k 1_{1\times n} K^P 1_{n\times 1}
$, which implies $k=-(\sum_{i=1}^n I^{\text{inj}})/(\sum_{i=1}^n K^P_i)$, concluding the proof. \hfill  $\qed$
\end{pf}

%%%%%%%%%%%%%%%%%%%%%%%%%%%%%%%%%%%%%%%%
%%%%%%%%%%%%%%%%%%%%%%%%%%%%%%%%%%%%%%%%
\section{Simulations}
\label{sec:simulations}
\begin{figure}
\begin{center}
\includegraphics[width=0.45\textwidth]{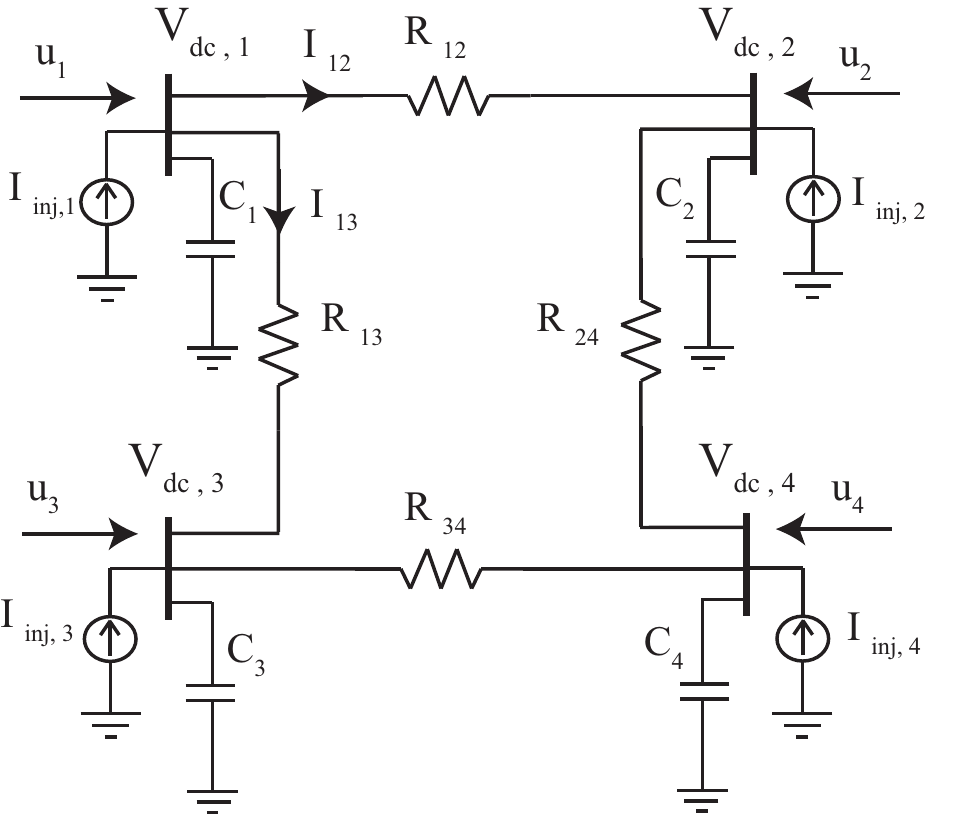}
\end{center}
\caption{Model and topology of the MTDC system considered in the simulations.}
\label{fig:system}
\end{figure}

Simulations of an MTDC system were conducted using MATLAB. The MTDC  was modelled by \eqref{eq:hvdc_scalar}, with $u_i$ given by either the droop controller \eqref{eq:droop}, or the distributed controller \eqref{eq:distributed_voltage_control}. The nominal voltage is assumed to be given by $V^{\text{nom}}= 100$ kV. The topology of the MTDC system is given by Figure \ref{fig:system}. The capacities are assumed to be $C_i=123.79 \; \mu \text{F}$ for $i=1,2,3,4$, while the resistances are assumed to be $R_{12}= 0.0154 \; \Omega$, $R_{13}= 0.0015 \; \Omega$, $R_{24}= 0.0015 \; \Omega$ and $R_{34}= 0.0154 \; \Omega$. The gains were set to $K^P_i= 10 \; \Omega^{-1}$ for $i=1,2,3,4$, and for both the VDM controller and the distributed controller. The remaining controller parameters were set to $\gamma = 0.005$ and $c_{ij} = R_{ij}^{-1} \; \Omega^{-1}$ for all $(i,j)\in \mathcal{E}$. Due to the long geographical distances between the DC converters, communication between neighboring nodes is assumed to be delayed with delay $\tau$ for the distributed controller. While the nominal system without time-delays is verified to be stable according to Theorem \ref{th:distributed_voltage_control},  time-delays might destabilize the system. It is thus of importance to study the effects of time-delays further. 
 The dynamics of the system \eqref{eq:hvdc_scalar} with the controller \eqref{eq:distributed_voltage_control}, and with time delay $\tau$ thus become

\begin{figure*}[ht]
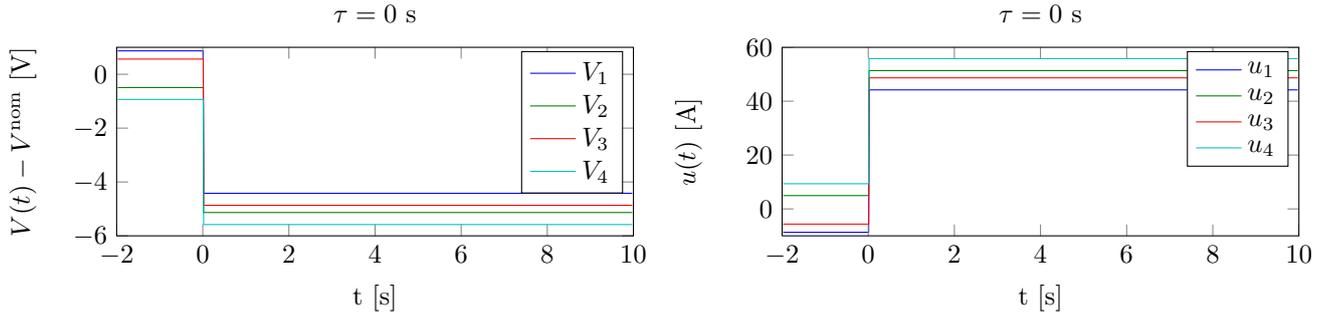

	\centering
$
\begin{array}{cc}
	\input{Simulations/V_tau_droop=0gamma=0005.tikz} & \input{Simulations/u_tau_droop=0gamma=0005.tikz}
\end{array}
$
\caption{The figure shows the voltages relative to the nominal voltage ($V_i-V^{\text{nom}}$), and the controlled injected currents $u_i$. The system model is given by \eqref{eq:hvdc_scalar}, and $u_i$ is given by the VDM controller \eqref{eq:droop}. The voltages and injected currents converge quickly to their stationary values. However, all voltages are below the nominal voltage, and the controlled injected currents are not equal.}
\label{fig:powersystems_sim_1_droop}
\end{figure*}

\begin{figure*}[ht]
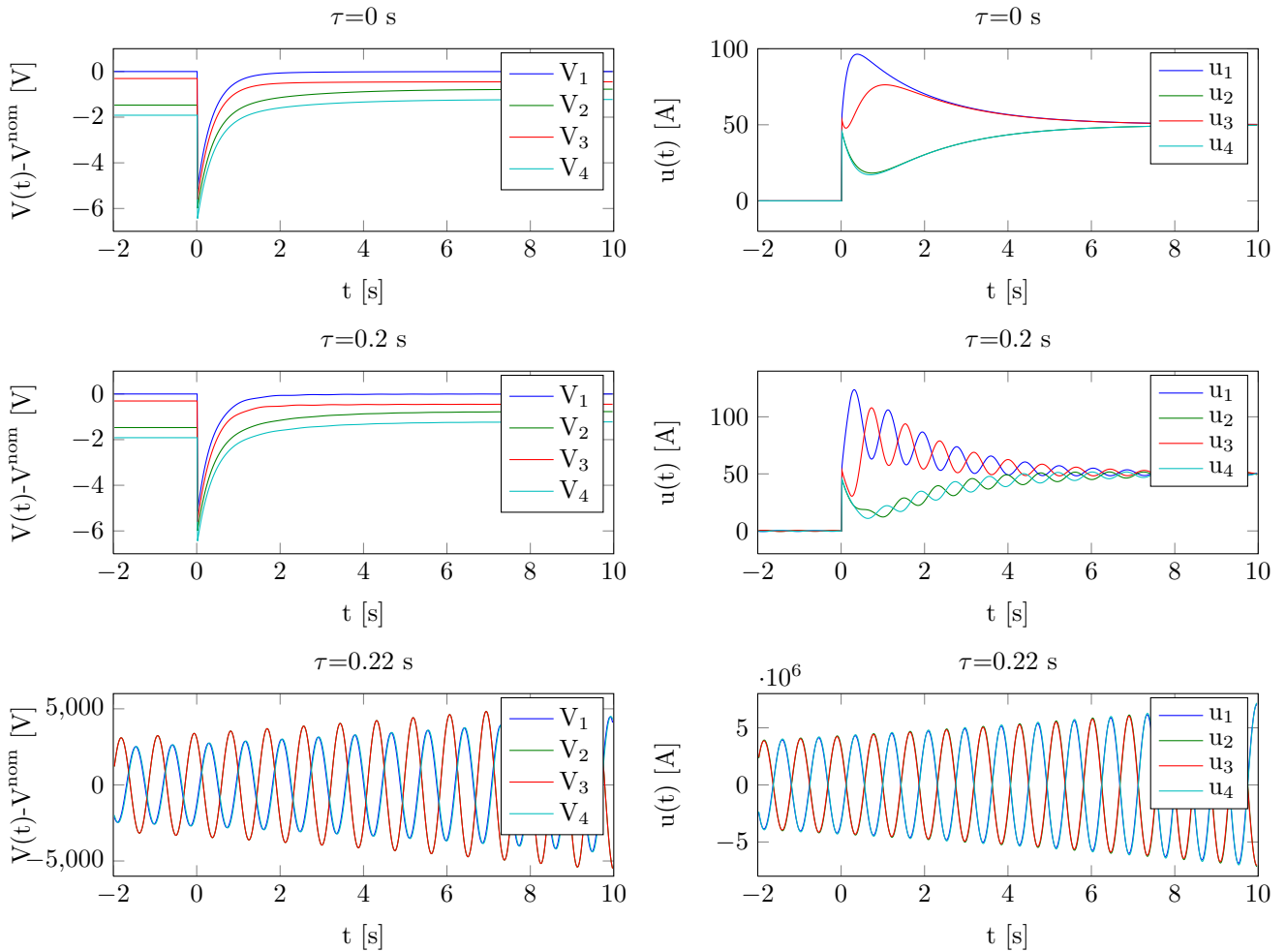

	\centering
$
\begin{array}{cc}
	\input{Simulations/V_tau=0gamma=0005.tikz} & \input{Simulations/u_tau=0gamma=0005.tikz} \\
	\input{Simulations/V_tau=02gamma=0005.tikz} & \input{Simulations/u_tau=02gamma=0005.tikz} \\
	\input{Simulations/V_tau=03gamma=0005.tikz} & \input{Simulations/u_tau=03gamma=0005.tikz}
\end{array}
$
\caption{The figure shows the voltages relative to the nominal voltage ($V_i-V^{\text{nom}}$), and the controlled injected currents $u_i$ of the converters for different time-delays $\tau$ on the communication links. The system model is given by \eqref{eq:hvdc_scalar}, and $u_i$ is given by the distributed controller \eqref{eq:distributed_voltage_control_delay}. The convergence times of the voltages and injected currents are in the order of a few seconds. On the other hand, we see that $V_1$ converges to $V^{\text{nom}}$, and that the controlled injected currents all converge to the same value.}
\label{fig:powersystems_sim_1}
\end{figure*}

\begin{align}
\label{eq:distributed_voltage_control_delay}
\begin{aligned}
u_i &= K^P(\hat{V}_i(t) -V_i(t)) \\
\dot{\hat{V}}_i &= K^V_i(V^{\text{nom}}-V_i(t)) \\
- &\gamma \sum_{j\in \mathcal{N}_i} c_{ij} \left( (\hat{V}_i(t') -V_i(t')){-}(\hat{V}_j(t'){-}V_j(t')) \right),
\end{aligned}
\end{align}
where $t'=t-\tau$. The injected currents are assumed to be initially given by $I^{\text{inj}}=[300,200,-100,-400]^T$ A, and the system is allowed to converge to the stationary solution. Since the injected currents satisfy $I_i^{\text{inj}}=0$, $u_i=0$ for $i=1,2,3,4$ by Theorem \ref{th:distributed_voltage_control}.
Then, at time $t=0$, the injected currents are changed due to changed power loads. The new injected currents are given by $I^{\text{inj}}=[  300, 200, -300, -400]^T$ A, i.e., only the injected current of converter $3$ is changed. The step response of the voltages $V_i$ and the controlled injected currents $u_i$ are shown in Figure \ref{fig:powersystems_sim_1_droop} for the droop controller \eqref{eq:droop}, and in Figure \ref{fig:powersystems_sim_1} for the distributed controller with time-delays \eqref{eq:distributed_voltage_control_delay}.

For the droop controller \eqref{eq:droop}, the system is stable as shown in Theorem \ref{th:droop_stability}. However, none of the voltages converges to  $V^{\text{nom}}$, and the controlled injected currents for the different converters do not converge to the same value, in accordance with Theorem \ref{th:droop_stability}. 

For the distributed controller \eqref{eq:distributed_voltage_control_delay} without delays, i.e., $\tau=0$ s, the voltages $V_i$ are restored to their new stationary values within $2$ seconds. The controlled injected currents $u_i$ converge to their stationary values within $8$ seconds. The simulations with time delays $\tau=0.1,0.22$ s, show that the controller is robust to moderate time-delays, but eventually the closed loop system becomes unstable. 

\setlength\fheight{2.5cm}
\setlength\fwidth{6.8cm}

\setlength\fheight{2.5cm}
\setlength\fwidth{6.8cm}

%%%%%%%%%%%%%%%%%%%%%%%%%%%%%%%%%%%%%%%%
%%%%%%%%%%%%%%%%%%%%%%%%%%%%%%%%%%%%%%%%
\section{Discussion and Conclusions}
\label{sec:discussion}
In this paper we have studied control of MTDC systems.
We have showed that a simple droop controller cannot satisfy the control objectives of voltage regulation and power sharing simultaneously, i.e., the controlled injected currents having a predefined ratio.  
We have proposed a distributed voltage controller for MTDC networks. We show that under mild conditions, there always exist controller parameters such that the closed-loop system is stable. In contrast to a decentralized droop controller, the proposed distributed controller is able to maintain the voltage levels of the converters close to the nominal voltages, while the injected current is shared proportionally amongst the converters.
We have validated our results through simulations, further showing that the distributed controller is robust to moderate time-delays. Future work will focus on finding upper bounds for the time-delay, guaranteeing closed loop stability under the distributed controller. 
\bibliography{references}
\bibliographystyle{plain}
\end{document}